\documentclass[12pt,leqno,twoside]{article} 
\usepackage{bezier,eufrak}
\input{rogg_kue.sty}
\begin{document}
\title{A relative Yoneda Lemma (manuscript)}
\author{M.\ K\"unzer}
\maketitle

\vspace{2cm}

\begin{abstract}
We construct set-valued right Kan extensions via a relative Yoneda Lemma.
\end{abstract}

\setcounter{section}{-1}

\section*{A remark of the referee}

As the referee pointed out, (\ref{ThK1}) `can essentially be found in much greater generality' in

\begin{itemize}
\item[{[Ke 82]}] {\sc G.\ M.\ Kelly,} {\it Basic concepts of enriched category theory,} LMS Lecture Notes 64, Cambridge University Press, 1982.
\end{itemize}

He continued to explain that to this end, one reformulates the formula [Ke 82, 4.6 (ii), p.\ 113], given in terms of weighted limits, by means of [Ke 82, 3.10, p.\ 99] and [Ke 82, 2.2, p.\ 58].

Therefore, we withdraw this note as a preprint. Since (\ref{ThK1}, \ref{ThK3}) might be of some use for the working mathematician, we leave it accessible as a manuscript. In using (\ref{ThK1}), 
the reader should refer to [Ke 82], in using (\ref{ThK3}), he should refer to [K 58].

\section{Introduction}

\subsection{A relative Yoneda Lemma}

The category of set-valued presheaves on a category $C$ shall be denoted by $C^\wdg$. The Yoneda embedding, sending $x$ to $\liu{C}{(-,x)}$, shall be denoted by $C\lraa{y_C} C^\wdg$. The presheaf 
category construction being contravariantly functorial, we obtain the

{\bf Proposition} (\ref{ThK1}, the relative Yoneda Lemma). {\it Given a functor $C\lraa{f} D$, we have
\[
f^\wedge\adj y_D^\wdg\c f^{\wdg\wdg}\c y_{C^\wdg}.
\]
There is a set theoretical caveat. In particular, this formula is correct only after some additional comments, see (\ref{ThK1}).} 

The right hand side is also known as the set-valued right Kan extension functor along $C^o\lraa{f^o} D^o$. There are several formulas for right Kan extensions in the literature, for instance 
using ends [ML 71, X.4], all of them necessarily yielding the same result up to natural equivalence, by uniqueness of the adjoint. In particular, (\ref{ThK1}) is merely still another
such a formula.

Letting $f = 1_C$, we recover the (absolute) Yoneda Lemma, thus giving a solution to the exercise [ML 71, X.7, ex.\ 2].
Concerning its origin, {\sc Mac Lane} recalled the following incident, taking place in around 1954/55 [ML 98].

\begin{quote}
{\sc Mac Lane,} then visiting Paris, was anxious to learn from {\sc Yoneda,} and commenced an interview with {\sc Yoneda} in a caf\'e at the Gare du Nord. The interview was continued on 
{\sc Yoneda}'s train until its departure. In its course, {\sc Mac Lane} learned about the lemma and subsequently baptized it.
\end{quote}

\subsection{A relative co-Yoneda Lemma}

Suppose given a functor $C\lraa{k}Z$.
{\sc Kan} constructed in [K 58, Th.\ 14.1] the left adjoint to the functor from $Z$ to $C^\wdg$ that sends $z\in Z$ to $\liu{Z}(k(-),z)\in C^\wdg$.
Now given a functor $C\lraa{f} D$, we can specialize to $Z = D^\wdg$ and to a functor $k$ that sends $c\in C$ to $\liu{D}{(-,fc)}\in D^\wdg$, thus obtaining a left adjoint to 
$C^\wdg\llaa{f^\wdg} D^\wdg$. ({\sc Kan}'s notation is as follows. Identify $\Zl = Z$, $\Vl = C$, $\Mfk^\Vl = C^\wdg$. The functor $H^\Vl(\Zl_\Vl,\Zl)$ maps from $\Zl_\Vl,\Zl$ to $\Mfk^\Vl$, 
i.e.\ from (the category of covariant functors from $\Vl$ to $\Zl$) $\ti\; \Zl$ to the category of presheaves on $\Vl$. Given $\Vl\lraa{k}\Zl$ and an object $z\in\Zl$, the pair $(k,z)$ is mapped 
to $\liu{\Zl}{(-,z)}\c k$.) We shall rephrase this special case of {\sc Kan}'s formula as follows, in order to be able to compare the left and the right adjoint of $f^\wdg$.

Let $C^\vee := (C^o)^\wdg$. The co-Yoneda embedding, sending $x$ to $\liu{C}{(x,-)}$, shall be denoted by $C\lraa{y'_C} C^\vee$. There is a tensor product $C^\wdg\ti C^\vee\lra (\Set)$, sending 
$v\ti v'$ to $v\ts_C v'$. The according univalent functor that sends $v$ to $v\ts_C -$ shall be denoted by $C^\wdg\lraa{z_C} C^{\vee\vee}$.

{\bf Proposition} ([K 58], cf.\ \ref{ThK3}, the relative co-Yoneda Lemma). {\it Given a functor $C\lraa{f} D$, we have
\[
y'_D\!^\vee\c f^{\vee\vee}\c z_C\adj f^\wdg.
\]
Again, there is a set-theoretical caveat.}

\subsection{Acknowledgements}

I would like to thank {\sc F.\ Borceux} for pointing out set theoretical difficulties. I would like to thank {\sc R.\ Rentschler} for kind hospitality. The author received support from the EU 
TMR-network `Algebraic Lie Representations', grant no.\ ERB FMRX-CT97-0100. 

\subsection{Notation}

\begin{footnotesize}
\begin{itemize}
\item[(i)] We write composition of morphisms on the right, $\lraa{a}\lraa{b} = \lraa{ab}$. However, we write composition of functors on the left, 
$\lraa{f}\lraa{g} = \lraa{g\c f}$. 
\item[(ii)] Given a functor $C\lraa{f} D$, the opposite functor between the opposite categories is denoted by $C^o\lraa{f^o}D^o$. 
\item[(iii)] Given functors $C\lradoublea{f}{f'} D\lradoublea{g}{g'} E$\vspace{3mm} and natural transformations $f\lraa{\alpha} f'$ and $g\lraa{\beta} g'$, we denote by 
$g\c f\lraa{\beta\c f} g'\c f$ the natural transformation defined by $(\beta\c f)c = \beta (fc)$ for $c\in C$ and by $g\c f\lraa{g\c\alpha} g\c f'$ the natural transformation defined by 
$(g\c\alpha)c = g(\alpha c)$ for $c\in C$. More generally,
\[
(g\c f\lraa{\beta\c\alpha} g'\c f') := (g\c f\lraa{\beta\c f} g'\c f\lraa{g'\c\alpha} g'\c f') = (g\c f\lraa{g\c\alpha} g\c f'\lraa{\beta\c f'} g'\c f').
\]
\end{itemize}
\end{footnotesize}          
\section{Universes}
\label{SecUniv}

\begin{quote}
\begin{footnotesize}
Since we will iterate the construction `forming the presheaf category over a category' once, we shall work in the setting of universes, which enables us to do such `large' constructions
when keeping track of the universe needed. Therefore, we start with a preliminary section to recall this well-known technique from [SGA 4 I, App.] and to fix some notation. 
\end{footnotesize}
\end{quote}

\begin{Definition}[{{\sc N.\ Bourbaki,} [SGA 4 I, App., 1, D\'ef.\ 1]}]
\label{DefU1}
\Absatz
A {\rm universe} is a set $\Ufk$ that satisfies the conditions (U1-4).
\begin{itemize}
\item[(U1)] $x\in y\in\Ufk$ implies $x\in\Ufk$.
\item[(U2)] $x,y\in\Ufk$ implies $\{x,y\}\in\Ufk$.
\item[(U3)] $x\in\Ufk$ implies $\Pfk(x)\in\Ufk$.
\item[(U4)] Given $I\in\Ufk$ and a map $I\lra\Ufk$, $i\lramaps x_i$, the union $\Cup_{i\in I} x_i$ is in $\Ufk$.
\end{itemize}
Here $\Pfk(x)$ denotes the power set of $x$.
\end{Definition}

\begin{Remark}
\label{RemU2}
\Absit\rm
\begin{itemize}
\item[(i)] $x\in\Ufk$ implies $\{ x\} = \{x,x\}\in\Ufk$.
\item[(ii)] $x\tm y\in\Ufk$ implies $x\in\Ufk$. In particular, if $y$ surjects onto some set $z$, then $z$ is in bijection to an element of $\Ufk$.
\item[(iii)] Let $(x,y) := \{x,\{x,y\}\}$. The element $x$ is the unique element of 
\[
\{ a\in (x,y) \; |\; \mb{for all }\; b\in (x,y)\ohne\{a\} \mb{ we have } a\in b\}, 
\]
since $\{x,y\}\in x$ would contradict {\sc von Neumann}'s axiom, asserting that any nonempty set $S$ contains an element that has an empty intersection with $S$. Now $\{x,y\}$ is the unique 
element of $(x,y)\ohne \{x\}$, and $y$ is the unique element of $\{x,y\}\ohne \{ x\}$.
\item[(iv)] $X\in\Ufk$ and $Y\in\Ufk$ implies $X\ti Y:= \Cup_{x\in X}\Cup_{y\in Y} \{ (x,y)\} \in\Ufk$.
\item[(v)] Given $I\in\Ufk$ and a map $I\lra\Ufk$, $i\lramaps x_i$, the disjoint union $\coprod_{i\in I} x_i := \Cup_{i\in I} X_i\ti\{ i\}$ is in $\Ufk$.
\item[(vi)] Given $I\in\Ufk$ and a map $I\lra\Ufk$, $i\lramaps x_i$, the product $\prod_{i\in I} x_i \tm \Pfk(\coprod_{i\in I} x_i)$ is in $\Ufk$.
\end{itemize}
\end{Remark}

Assume given universes $\Ufk$, $\Vfk$ and $\Wfk$ such that 
\[
\leer\neq\Ufk\in\Vfk\in\Wfk. 
\]
This may be achieved by means of {\sc Bourbaki}'s axiom (A.6) in [SGA 4 I, App., 4], which says that any set is element of some universe. Note that $\Ufk\tm\Vfk\tm\Wfk$.

For set theoretical purposes the following definition is convenient. 

\begin{Definition}[{cf.\ [M 65, I.2]}]
\label{DefU3}
A {\rm small category} $C$ is a tuple $(M,K\tm M\ti M\ti M)$, $M$ being a set, subject to the conditions (i-iii). We write $M =: \mb{\rm Mor}\, C$, the set of {\rm morphisms} of $C$. $K$ is the 
{\rm composition law} of $C$.

\begin{itemize}
\item[(i)] For every $(a,b)\in M\ti M$ there exists at most one $c\in M$ such that $(a,b,c)\in K$. In this case, we write $c = ab$ and say that `(the {\rm composite}) $ab$ exists'. 
\end{itemize}
Let 
\[
\begin{array}{rcl}
\Ob C := \{\; a\in M & | & \mb{for all $b$ for which $ba$ exists, we have $ba = b$, and}\\
                     &   & \mb{for all $c$ for which $ac$ exists, we have $ac = c$.} \}\;\tm\; \mb{\rm Mor}\, C \\
\end{array}
\]
be the set of {\rm objects} of $C$. If $a\in\Ob C$, we also write $a =: 1_a$. 
\begin{itemize}
\item[(ii)] For every $a\in M$ there exists a unique {\rm source} element $s(a)\in\Ob C$ such that $s(a) a$ exists, and a unique {\rm target} element $t(a)\in\Ob C$ such that $a t(a)$ exists.

\item[(iii)] For $a,b\in M$ such that the composite $ab$ exists, the composites $t(a) b$ and $a s(b)$ exist.

\item[(iv)] Given $a,b,c\in M$ such that $ab$ and $bc$ exist. Then $(ab)c$ and $a(bc)$ exist and equal each other. 
\end{itemize}

The set of morphisms with start $x$ and target $y$ shall be denoted by $\liu{C}{(x,y)}$.

\rm
Suppose given $a,b\in M$. Note that if $ab$ exists, we have $s(ab) = s(a)$ and $t(ab) = t(b)$. Note that $s(s(a)) = s(a)$, $t(t(a)) = t(a)$. Note that $ab$ exists iff $t(a) = s(b)$.
\end{Definition}

\begin{Definition}
\label{DefU3_5}
A set $X$ is said to be {\rm $\Ufk$-small} if there exists a bijection from $X$ to an element of $\Ufk$. 
The category of $\Ufk$-small sets is denoted by $(\Set_\Ufk)$. Note that by a skeleton argument, $(\Set_\Ufk)$ is equivalent to the category of sets contained as elements in $\Ufk$, denoted by 
$(\Set^0_\Ufk)$.

A small category $C$ is said to be {\rm $\Ufk$-small} if $\mb{\rm Mor}\, C$ is $\Ufk$-small. A category $C$ is said to be {\rm essentially $\Ufk$-small} if it is equivalent to a $\Ufk$-small 
category, or, equivalently, if it has a $\Ufk$-small skeleton. 
\end{Definition}

\begin{Remark}
\label{RemU4}
\rm
The category $(\Set^0_\Ufk)$ is $\Vfk$-small since $\coprodd{(x,y)\in\Ufk\ti\Ufk} \liu{(\sSet^0_\Ufk)}{(x,y)}$ is an element of $\Vfk$. Hence the category $(\Set_\Ufk)$ is essentially $\Vfk$-small.
\end{Remark}

Given categories $C$ and $D$, we denote the category of functors mapping from $C$ to $D$ by $\bo C, D\bc$.

\begin{Lemma}
\label{RemU5}
If $C$ and $D$ are $\Ufk$-small, then $\bo C,D\bc$ is $\Ufk$-small. If $C$ and $D$ are essentially $\Ufk$-small, then $\bo C,D\bc$ is essentially $\Ufk$-small.

\rm
Using induced equivalences, it suffices to prove the first assertion. But then 
\[
\Ob\bo C,D\bc \tm \Pfk(\mb{Mor}\, C\ti\mb{Mor}\, D)
\]
is $\Ufk$-small. Moreover, given $f,g\in\bo C,D\bc$, the set 
\[
\liu{\sbo C,D\sbc}{(f,g)}\tm\coprodd{x\in\sOb C}\liu{D}{(fx,gx)}
\]
is $\Ufk$-small. Hence $\mb{Mor}\,\bo C,D\bc$ is $\Ufk$-small.
\end{Lemma}

\begin{Lemma}
\label{LemU6}
Given an essentially $\Ufk$-small category $C$, the category 
\[
C^{\wdg_\Ufk} := \bo C^o,(\Set_\Ufk)\bc
\]
of {\rm $(\Set_\Ufk)$-valued presheaves} over $C$ is essentially $\Vfk$-small. Sometimes we write $C^\wdg = C^{\wdg_\Ufk}$ if the universe is unambiguous. Likewise, the category 
\[
C^{\vee_\Ufk} := \bo C,(\Set_\Ufk)\bc
\]
of {\rm $(\Set_\Ufk)$-valued copresheaves} over $C$ is essentially $\Vfk$-small. Sometimes we write $C^\vee = C^{\vee_\Ufk}$. 

\rm
This follows from (\ref{RemU4}, \ref{RemU5}). 
\end{Lemma}

\begin{Definition}
Given an essentially $\Ufk$-small category $C$, we have the {\rm Yoneda embedding}
\[
\begin{array}{rcl}
C & \lraa{y_C} & C^{\wdg_\Ufk} \\
x & \lra       & \liu{C}{(-,x)} \\
\end{array}
\]
and the {\rm co-Yoneda embedding}
\[
\begin{array}{rcl}
C & \lraa{y'_C} & C^{\vee_\Ufk} \\
x & \lra        & \liu{C}{(x,-)}. \\
\end{array}
\]

For a functor $C\lraa{f} D$, we denote 
\[
\begin{array}{rcl}
C^{\wdg_\Ufk}                        & \llaa{f^\wdg} & D^{\wdg_\Ufk} \\
(u\c f^o\lraa{\beta\c f^o} u'\c f^o) & \lla          & (u\lraa{\beta} u'). \\
\end{array}
\]
Given functors $C\lradoublea{f}{g} D$\vspace{2mm} and a natural transformation $f\lraa{\alpha} g$, we denote by $f^\wdg\llaa{\alpha^\wdg} g^\wdg$ the natural transformation that is given at 
$u\in D^{\wdg_\Ufk}$ by the morphism $u\c f^o \llaa{u\c \alpha^o} u\c g^o$ in $C^{\wdg_\Ufk}$, that evaluated at $c\in C^o$ in turn yields $ufc\llaa{u\alpha c} ugc$. 

Analogously for $C^{\vee_\Ufk}\llaa{f^\vee}D^{\vee_\Ufk}$ and $f^\vee\lraa{\alpha^\vee} g^\vee$.
\end{Definition}

\section{The right Kan extension}
\label{SecRight}

\begin{Proposition}[the relative Yoneda Lemma]
\label{ThK1}
Given essentially $\Ufk$-small categories $C$ and $D$, and a functor $C\lraa{f} D$. Then the right adjoint $C^{\wdg_\Ufk}\lraa{\phi} D^{\wdg_\Ufk}$ of 
$C^{\wdg_\Ufk}\llaa{f^\wdg}D^{\wdg_\Ufk}$ is given by 
\begin{center}
\begin{picture}(250,450)
\put( -80, 400){$C^{\wdg_\Ufk\wdg_\Vfk}$}
\put(  50, 410){\vector(1,0){130}}
\put(  90, 425){$\scm f^{\wdg\wdg}$}
\put( 200, 400){$D^{\wdg_\Ufk\wdg_\Vfk}$}
\put( -10,  50){\vector(0,1){330}}
\put( -70, 200){$\scm y_{C^\wdg}$}
\put( 220,  60){\vector(0,1){120}}
\put( 210,  60){\oval(20,20)[b]}
\put( 200, 200){$D^{\wdg_\Vfk}$}
\put( 220, 380){\vector(0,-1){130}}
\put( 230, 310){$\scm y_D^\wdg$}
\put( -40,   0){$C^{\wdg_\Ufk}$}
\put(  50,  10){\vector(1,0){130}}
\put( 100,  25){$\scm \phi$}
\put( 200,   0){$D^{\wdg_\Ufk}$}
\end{picture}
\end{center}
Keeping the name of the functor after restricting the image to $D^{\wdg_\Ufk}$, we write shorthand
\[
\fbox{$\;
f^\wdg\adj y_D^\wdg\c f^{\wdg\wdg}\c y_{C^\wdg}.
\;$}
\]
The unit of this adjunction
\[
\fbox{$\;
1_{D^\wdg}\lraa{\eps} \phi\c f^\wdg
\;$}
\]
at $u\in D^{\wdg_\Ufk}$, i.e.\ 
\[
u\lraa{\eps u} \liu{C^\wdg}{(f^{\wdg,o}\circ y_D^o(-),u\c f^o)},
\]
applied to $d\in D$, is given by
\[
\begin{array}{rcl}
ud & \lraa{\eps u d} & \liu{C^\wdg}{(\liu{D}{(f^o(-),d)},u\c f^o)} \\
x  & \lra            & (x)\eps u d, \\
\end{array}
\]
where the natural transformation $(x)\eps u d$ sends at $c\in C$ 
\[
\begin{array}{rcl}
\liu{D}{(fc,d)} & \lraa{(x)\eps u d} & u f c \\
a               & \lra               & (a)\left[\Big((x)\eps u d\Big)c\right] = (x)(ua^o). \\
\end{array}
\]

The counit of this adjunction 
\[
\fbox{$\;
f^\wdg\c \phi \lraa{\eta} 1_{C^\wdg}
\;$}
\]
at $v\in C^\wdg$, i.e.\
\[
\liu{C^\wdg}{(f^{\wdg,o}\c y_D^o\c f^o(-),v)} \lraa{\eta v} v,
\]
applied to $c\in C$, is given by
\[
\begin{array}{rcl}
\liu{C^\wdg}{(\liu{D}(f^o(-),f c),v)} & \lraa{\eta v c} & vc \\
\xi                                   & \lra            & (\xi)\eta v c = (1_{f c})\xi c\\
\end{array}
\]

\rm
Since the set
\[
\liu{C^\wdg}{(\liu{D}{(f(-),d)},v)} \tm \prodd{c\in\sOb C} \liu{(\sSet_\Ufk)}{(\liu{D}{(fc,d)},vc)}
\]
is $\Ufk$-small, $\phi$ exists as the factorization of $y_D^\wdg\c f^{\wdg\wdg}\c y_{C^\wdg}$ over the inclusion $D^{\wdg_\Ufk}\tm D^{\wdg_\Vfk}$.

Various compatibilities need to be verified to ensure the well-definedness of $\eps$ and $\eta$ (\footnote{
Given $u\in D^{\wdg_\Ufk}$, we need to see that $u\lraa{\eps u}\phi\c f^\wdg(u)$ is a natural transformation. Suppose given $d'\lraa{a} d$ in $D$. We have to show that for any
$x\in ud$
\[
(x)\Big(\eps u d\Big)\Big((\phi\c f^\wdg)ua^o\Big) = (x)\Big(ua^o\Big)\Big(\eps ud'\Big)
\]
is an equality of natural transformations from $\liu{D}{(f^o(-),d')}$ to $u\c f^o$. At $c\in C$, the element $f c\lraa{a'} d'$ is mapped by the left hand side to $(x)(u(a'a)^o)$, and by the 
right hand side to $(x (ua^o))(ua'^o)$.

We need to see that $1_{D^{\wdg_\Ufk}}\lraa{\eps}\phi\c f^\wdg$ is a natural transformation. Suppose given $u\lraa{s} u'$ in $D^{\wdg_\Ufk}$. We have to show that for any $d\in D$ and any 
$x\in ud$
\[
(x) \Big(sd\Big) \Big(\eps u' d\Big) = (x)\Big(\eps u d\Big)\Big((\phi\c f^\wdg (s))d\Big)
\]
is an equality of natural transformations from $\liu{D}{(f^o(-),d)}$ to $u'\c f$. At $c\in C$, the element $f c\lraa{a} d$ is mapped by the left hand side to
$(x)(sd)(u'a^o)$, and by the right hand side to $(x)(u a^o)(s(f c))$.

Given $v\in C^{\wdg_\Ufk}$, we need to see that $f^\wdg\c\phi(v)\lraa{\eta v} v$ is a natural transformation. Suppose given $c'\lraa{b} c$ in $C$. We have to show that for any
$\xi\in (f^\wdg\c \phi(v))c = \liu{C^\wdg}{(\liu{D}(f^o(-),f c),v)}$ we have
\[
(x)\Big(\eta vc\Big)\Big(vb^o\Big) = (x)\Big((f^\wdg\c\phi(v))b^o\Big)\Big(\eta vc'\Big).
\]
The left hand side yields $(1_{fc})(\xi c)(v b^o)$, the right hand side yields $(fb)\xi c'$.

We need to see that $f^\wdg\c\phi\lraa{\eta} 1_{C^{\wdg_\Ufk}}$ is a natural transformation. Suppose given $v\lraa{t} v'$ in $C^{\wdg_\Ufk}$. We have to show that for any $c\in C$ and any 
$\xi\in (f^\wdg\c \phi(v))c = \liu{C^\wdg}{(\liu{D}(f^o(-),f c),v)}$ we have
\[
(\xi)\Big(\eta vc\Big)\Big(tc\Big) = (\xi)\Big((f^\wdg\c \phi(t))c\Big)\Big(\eta v'c\Big).
\]
The left hand side yields $(1_{f c})(\xi c)(tc)$. The right hand side yields $(1_{f c})((\xi t)c)$.
}).

We have to show that $(f^\wdg\c \eps)(\eta\c f^\wdg) = 1_{f^\wdg}$. Suppose given $u\in D^\wdg$, $c\in C$ and $x\in (f^\wdg u) c = u(f c)$. We obtain
\[
\begin{array}{rcl}
(x)\Big( (f^\wdg\c\eps)uc\Big)\Big((\eta\c f^\wdg)uc\Big)
& = & (x) \Big( \eps u (fc)\Big)\Big( \eta(u\c f^o)c\Big) \\
& = & (1_{fc})\left[\Big((x)(\eps u(fc))\Big) c\right] \\
& = & (x)(u 1_{fc}^o) \\
& = & x. \\
\end{array}
\]

We have to show that $(\eps\c\phi)(\phi\c\eta) = 1_{\phi}$. Suppose given $v\in C^{\wdg_\Ufk}$ and $d\in D$. Note that $\phi v = \liu{C^\wdg}{(f^{\wdg,o}\c y_D^o (-),v)}$ and therefore
$\phi v d = \liu{C^\wdg}{(\liu{D}{(f^o(-),d)},v)}$. The application $\Big((\eps\c\phi)vd\Big)\Big((\phi\c\eta)vd\Big)$ writes
\[
\liu{C}{(\liu{D}{(f^o(-),d)},v)}\lraa{\eps (\phi v)d} \liu{C^\wdg}{\Big(\liu{D}{(f^o(-),d),\liu{C^\wdg}{(f^{\wdg,o}\c y_D^o\c f^o(-),v)}}\Big)} 
\lraa{(-)\eta v} \liu{C^{\wdg}}{(\liu{D}{(f^o(-),d)},v)}
\]
An element $\xi\in\phi v d$, i.e.\ $\liu{D}{(f^o(-),d)}\lraa{\xi} v$, is thus mapped to the composite
\[
\liu{D}{(f^o(-),d)}\mra{(\xi)\eps (\phi v)d} \liu{C^\wdg}{(f^{\wdg,o}\c y_D^o \c f^o (-),v)} \lraa{\eta v} v.
\]
Now suppose given $c\in C$ and $fc\lraa{a}d$. We have
\[
\begin{array}{rclcl}
\liu{D}{(fc,d)} & \mra{((\xi)\eps (\phi v) d)c} & \liu{C^\wdg}{\Big(\liu{D}{(f^o(-),fc)},v\Big)}       & \lraa{\eta vc} & vc \\
a               & \mra{}                        & (\xi)((\phi v)a^o)                                   &                &    \\
                &                               & = (\xi)\liu{C^\wdg}{\Big(\liu{D}{(f^o(-),a)},v\Big)} &                &    \\
                &                               & = \liu{D}{(f^o(-),a)}\xi                             & \lra           & (1_{fc})\liu{D}{(fc,a)}(\xi c) \\
                &                               &                                                      &                & = (a)(\xi c), \\
\end{array}
\]
whence $\xi\Big((\eps\c\phi)vd\Big)\Big((\phi\c\eta)vd\Big) = \xi$.
\end{Proposition}

\begin{Remark}[the absolute Yoneda Lemma]
\label{RemK2}
In case $f = 1_C$, we obtain $1 \adj y_C^\wdg \c y_{C^\wdg}$, which by uniqueness of the right adjoint yields the comparison isomorphism
\[
\begin{array}{rcl}
\liu{C^\wdg}{(\liu{C}(-,c),v)} & \lraisoa{\eta v c} & vc \\
\xi                              & \lra               & (1_c)\xi c,\\
\end{array}
\]
at $v\in C^{\wdg_\Ufk}$ and $c\in C$, with inverse given by $\eps vc$.
\end{Remark}

\begin{Corollary}
\label{CorK2_5}
If $f$ is full, then the counit $\eta$ of the adjunction $f^\wdg\adj y_D^\wdg\c f^{\wdg\wdg}\c y_{C^\wdg}$ is a monomorphism. If $f$ is full and faithful, then $\eta$ is an isomorphism.
\end{Corollary}          
\section{The left Kan extension}
\label{SecLeft}

\begin{quote}
\begin{footnotesize}
For sake of comparison to (\ref{ThK1}), we rephrase the pertinent case of {\sc Kan}'s formula in our setting. 
\end{footnotesize}
\end{quote}

Let $C$ be a $\Ufk$-small category. Let $v\in C^{\wdg_\Ufk}$, let $w\in C^{\vee_\Ufk}$. We define the set $v\ts_C w$ as the quotient of the disjoint union
\[
v\ti_C w := \coprodd{c\in C} vc\ti wc
\]
modulo the equivalence relation generated by the following relation $\sim_C$. The equivalence class of $(p,q)\in vc\ti wc$, $c\in C$, shall be denoted by $p\ts q$.

Given $(p,q)\in vc\ti wc$, $(p',q')\in vc'\ti wc'$, we say that $(p, q)\sim_C (p', q')$ if there exists a morphism $c\lraa{a} c'$ such that
\[
\begin{array}{rcl}
(p')v a^o & = & p \\
(q) wa    & = & q'. \\
\end{array}
\]
Thus the quotient map $v\ti_C w\lraa{\nu} v\ts_C w$ has the following universal property. Given a map $v\ti_C w\lraa{\nu'} X$ such that for any morphism $c\lraa{a}c'$, any $p'\in vc'$ and 
any $q\in wc$ we have 
\[
((p')va^o, q)\nu' = (p', (q)wa)\nu',
\]
there exists a unique map $v\ts_C w\lraa{\w \nu'} X$ such that $\nu' = \nu\w\nu'$.

In particular, given morphisms $v\lraa{m} v'$ and $w\lraa{n}w'$, we obtain a map $m\ts_C n$ that maps an element represented by $(p,q)\in vc\ti wc$, $c\in C$, as follows.
\[
\begin{array}{rllcrll}
v & \ts_C & w & \lraa{m\ts_C n} & v'    & \ts_C & w' \\
p & \ts   & q & \lra            & (p)mc & \ts   & (q)nc \\
\end{array}
\]
Thus the tensor product defines a functor $C^{\wdg_\Ufk}\ti C^{\vee_\Ufk} \mra{=\,\ts_C\, -} (\Set_\Ufk)$. We denote the univalent tensor product functor by
\[
\begin{array}{rcl}
C^{\wdg_\Ufk} & \lraa{z_C} & C^{\vee_\Ufk\vee_\Ufk} \\ 
v             & \lra       & v\ts_C -. \\
\end{array}
\]

\begin{Proposition}[{the relative co-Yoneda Lemma, {\sc Kan} [K 58, Th.\ 14.1]}]
\label{ThK3}
Given $\Ufk$-small categories $C$ and $D$, and a functor $C\lraa{f} D$. The left adjoint $C^{\wdg_\Ufk}\lraa{\psi}D^{\wdg_\Ufk}$ of $C^{\wdg_\Ufk}\llaa{f^\wdg}D^{\wdg_\Ufk}$ is given by
\begin{center}
\begin{picture}(250,250)
\put( -80, 200){$C^{\vee_\Ufk\vee_\Ufk}$}
\put(  50, 210){\vector(1,0){130}}
\put(  90, 225){$\scm f^{\vee\vee}$}
\put( 200, 200){$D^{\vee_\Ufk\vee_\Ufk}$}
\put(   0,  50){\vector(0,1){130}}
\put( -40, 110){$\scm z_C$}
\put( 220, 180){\vector(0,-1){130}}
\put( 230, 110){$\scm y'_D\!^\wdg$}
\put( -40,   0){$C^{\wdg_\Ufk}$}
\put(  50,  10){\vector(1,0){130}}
\put( 100,  25){$\scm\psi$}
\put( 200,   0){$(D^o)^{\vee_\Ufk} = D^{\wdg_\Ufk}$}
\end{picture}
\end{center}
For short, 
\[
\fbox{$\;
y'_D\!^\vee\c f^{\vee\vee}\c z_C\adj f^\wdg.
\;$}
\]
The unit of this adjunction
\[
\fbox{$\;
1_{C^\wdg}\lraa{\eps} f^\wdg\c \psi
\;$}
\]
at $v\in C^\wdg$, i.e.\ 
\[
v\lraa{\eps v} v\ts_C f^\vee \c y'_D\c f^o (-),
\]
applied to $c\in C$, is given by
\[
\begin{array}{rcl}
vc & \lraa{\eps vc} & v\ts_C \liu{D}{(fc,f(-))} \\
x  & \lra           & x\ts 1_{fc}\; . \\
\end{array}
\]

The counit of this adjunction 
\[
\fbox{$\;
\psi\c f^\wdg\lraa{\eta} 1_{D^\wdg}
\;$}
\]
at $u\in D^\wdg$, i.e.\
\[
u\c f^o\ts_C f^\vee\c y'_D(-) \lraa{\eta u} u,
\]
applied to $d\in D$, is given by
\[
\begin{array}{rllcl}
u\c f^o & \ts_C & \liu{D}{(d,f(-))} & \lraa{\eta ud} & ud \\
p       & \ts   & q                 & \lra           & (p)uq^o,  \\
\end{array}
\]
where $p\ts q$ is represented by $(p,q)\in ufc\ti \liu{D}{(d,fc)}$ for some $c\in C$.

\rm
Various compatibilities have to be verified to ensure the well-definedness of $\eps$ and $\eta$ (\footnote{
Given $v\in C^\wdg$, we need to see that $v\lraa{\eps v} f^\wdg\c \psi(v)$ is a natural transformation. Suppose given $c'\lraa{b}c$. We have to show that for any $x\in vc$
\[
(x)\Big(\eps v c\Big)\Big((f^\vee\c\psi(v))b^o\Big) = (x)\Big(vb^o\Big)\Big(\eps vc'\Big).
\]
The left hand side yields $x\ts fb$. The right hand side yields $x(vb^o)\ts 1_{fc'}$.

We need to see that $1_{C^\wdg}\lraa{\eps} f^\wdg\c \psi$ is a natural transformation. Suppose given $v\lraa{t}v'$ in $C^\wdg$. For any $c\in C$ and any $x\in vc$ we have
\[
(x)\Big(tc\Big)\Big(\eps v'c\Big) = (x)\Big(\eps vc\Big)\Big(tc\ts\liu{D}(fc,f(-))\Big) = (x)tc\ts 1_{fc}.
\]

Given $u\in D^\wdg$ and $d\in D$, we need to see that $\eta ud$ is a well-defined map. Suppose given $c,c'\in C$, $c\lraa{b} c'$ in $C$ and $p'\in ufc'$, $q\in \liu{D}{(d,fc)}$. Since
\[
\Big((p')uf^o b^o\Big)uq^o = \Big(p'\Big) u(q(fb))^o,
\]
the universal property applies.

Given $u\in D^\wdg$, we need to see that $\psi\c f^\wdg(u)\lraa{\eta u} u$ is a natural transformation. Suppose given $d'\lraa{a} d$. For all $c\in C$, $p\in ufc$ and $q\in\liu{D}{(d,fc)}$
we obtain
\[
(p\ts q)\Big(\eta ud\Big)\Big(ua^o\Big) = (p\ts q)\Big((\psi\c f^\wdg(u))a^o\Big)\Big(\eta ud'\Big) = (p)u(aq)^o.
\]

We need to see that $\psi\c f^\wdg\lraa{\eta} 1_{D^\wdg}$ is a natural transformation. Suppose given $u'\lraa{s} u$ in $D^\wdg$. We have to show that for all $d\in D$, $c\in C$, $p\in ufc$ and 
$q\in\liu{D}{(d,fc)}$, we have
\[
(p\ts q)\Big(\eta ud\Big)\Big(sd\Big) = (p\ts q)\Big((\psi\c f^\wdg(s))d\Big)\Big(\eta u'd\Big).
\]
The left hand side yields $(p)(u q^o)(sd)$. The right hand side yields $(p)(s(fc))(u' q^o)$.
}).

We have to show that $(\eps\c f^\wdg)(f^\wdg\c\eta) = 1_{f^\wdg}$. Suppose given $u\in D^\wdg$, $c\in C$ and $x\in (f^\wdg u)c$. We obtain
\[
\begin{array}{rcl}
(x)\Big((\eps\c f^\wdg)u c\Big)\Big((f^\wdg\c\eta)uc\Big) 
& = & (x)\Big(\eps(u\c f^o) c\Big)\Big(\eta u(fc)\Big) \\
& = & (x\ts 1_{fc})\Big(\eta u(fc)\Big) \\
& = & (x) u 1_{fc}^o. \\
& = & x. \\
\end{array}
\]

We have to show that $(\psi\c\eps)(\eta\c\psi) = 1_\psi$. Suppose given $v\in C^\wdg$, $d\in D$, $c\in C$, $s\in vc$, $t\in\liu{D}{(d,fc)}$, so that $s\ts t\in \psi v d = v\ts_C \liu{D}{(d,f(-))}$.
We obtain
\[
\begin{array}{rcl}
(s\ts t)\Big((\psi\c\eps)vd\Big)\Big((\eta\c\psi)vd\Big) 
& = & (s\ts t)\Big(\eps v\ts \liu{D}{(d,f(-))}\Big)\Big(\eta(v\ts_C f^\vee\c y'_D)d\Big) \\
& = & \Big((s\ts 1_{fc})\ts t\Big)\Big(\eta(v\ts_C f^\vee\c y'_D)d\Big) \\
& = & (s\ts 1_{fc}) \Big(v\ts_C \liu{D}{(t,f(-))}\Big) \\
& = & (s\ts t). \\
\end{array}
\] 
\end{Proposition}

\begin{Remark}[the absolute co-Yoneda Lemma]
\label{RemK4}
In case $f = 1_C$, we obtain $y'_D\!^\vee\c z_C\adj 1_{C^\wdg}$, which by uniqueness of the left adjoint yields the comparison isomorphism
\[
\begin{array}{rllcl}
v & \ts_C & \liu{C}{(c,-)} & \lraisoa{\eta vc} & vc       \\
s & \ts   & t              & \lra              & (s)vt^o  \\
\end{array}
\]
at $v\in C^{\wdg_\Ufk}$ and $c\in C$, with inverse given by $\eps v c$.
\end{Remark}

\begin{Corollary}
\label{CorK5}
If $f$ is full, then the unit $\eps$ of the adjunction $y'_D\!^\vee\c f^{\vee\vee}\c z_C\adj f^\wdg$ is an epimorphism. If $f$ is full and faithful, then $\eps$ is an isomorphism.
\end{Corollary}          
\section{References}

\begin{footnotesize}
{\sc Bourbaki, N.}\upl
\begin{itemize}
\item[{[SGA 4]}] {\it Univers,} appendice de l'expos\'e I de SGA 4, tome 1, SLN 269, 1972.
\end{itemize}

{\sc Kan, D.\ M.}\upl
\begin{itemize}
\item[{[K 58]}] {\it Adjoint functors,} Trans.\ Am.\ Math.\ Soc.\ 87, pp.\ 294-329, 1958.
\end{itemize}

{\sc Mac Lane, S.}\upl
\begin{itemize}
\item[{[ML 71]}] {\it Categories for the Working Mathematician,} GTM 5, 1971.
\item[{[ML 98]}] {\it The Yoneda Lemma,} Math.\ Japon.\ 47 (1), p.\ 156, 1998.
\end{itemize}

{\sc Mitchell, B.}\upl
\begin{itemize}
\item[{[M 65]}] {\it Theory of categories,} Academic Press, 1965.
\end{itemize}

\vspace*{8cm}

\begin{flushleft}
Matthias K\"unzer\\
Universit\"at Ulm\\
Abt. Reine Mathematik\\
D-89069 Ulm\\
Germany\\
kuenzer@mathematik.uni-ulm.de\\
\end{flushleft}
\end{footnotesize}
\end{document}